\documentclass{article}
\usepackage{amsmath}

\numberwithin{equation}{section}
\newtheorem{thm}{Theorem}[section]
\newtheorem{cor}[thm]{Corollary}
\newtheorem{prop}[thm]{Proposition}

\newtheorem{lem}[thm]{Lemma}

\newcommand{\be}{\begin{equation}}
\newcommand{\ee}{\end{equation}}
\newcommand{\ben}{\begin{enumerate}}
\newcommand{\een}{\end{enumerate}}
\newcommand{\beq}{\begin{eqnarray}}
\newcommand{\eeq}{\end{eqnarray}}
\newcommand{\beqn}{\begin{eqnarray*}}
\newcommand{\eeqn}{\end{eqnarray*}}

\newcommand{\Qed}{\hspace*{\fill}Q.E.D.}  

\title{\Large  On a class of projectively flat Finsler metrics}
\author{Benling Li\footnote{Research is supported by the NNSFC(11371209), ZPNSFC(LY13A010013), Scientific Research Fund of Zhejiang Provincial Education Department (Y201326630),  Ningbo Natural Science Foundation (2013A610101) and K.C. Wong Magna Fund in Ningbo University.} \    and Zhongmin Shen}

\date{November 10, 2014}
\begin{document}
\maketitle
\begin{abstract}
 In this paper,  we study a class of Finsler metrics composed by a Riemann metric $\alpha=\sqrt{a_{ij}(x)y^i y^j}$ and a $1$-form $\beta=b_i(x)y^i$ called general ($\alpha$, $\beta$)-metrics.
We classify those projectively flat when $\alpha$ is projectively flat. By solving the corresponding nonlinear PDEs, the metrics in this class are totally determined.
Then a new group of  projectively flat Finsler metrics is found.
\end{abstract}


\textbf{2000 MR Subject Classification.}  53B40, 53C60.

\section{Introduction}

   The regular case of the famous Hilbert's Fourth Problem is to characterize Finsler
metrics on an open subset in $R^n$ whose geodesics are straight lines as a point set.
Such Finsler metrics are called {\it projectively flat} Finsler metrics.
In the past years, many non-trivial(non-Minkowskian) projectively flat Finsler metrics have been found.
The simplest class is projectively flat Riemannian metrics. As we known, they are equivalent to
Riemannian metrics with constant sectional curvature by Beltrami's theorem.
However, it is not true in the non-Riemannian case.

Recently, in \cite{Li1} the first author studied  dually flat Finsler metrics arisen
from information geometry and found that many dually flat Finsler metrics
can be constructed by projectively flat Finsler metrics.
More information of dually flat Finsler metrics can be found in \cite{Cheng}\cite{Huang}\cite{Xia}.
These also motivate us to search more projectively flat Finsler metrics.
Many known examples are related to Riemannian metrics,
such as the famous Funk metric $ \Theta = \Theta(x,y) $ found in 1904 and the Berwald's metric $B = B(x,y)$ found in 1929.
This leads to the  study of {\it ($\alpha$, $\beta$)-metric}
defined by a Riemannian metric $\alpha$ and a $1$-form $\beta$,
\be
F =\alpha \phi(s), \ \ \ \ \ s=\frac{\beta}{\alpha},\label{Fber***}
\ee
where $\phi=\phi(s)$ satisfies  certain condition such that $F$ is a (positive definite) Finsler metric. The second author gave the equivalent conditions of projectively flat ($\alpha$, $\beta$)-metrics in \cite{ShenProj1}. Later on, we classified those with constant flag curvature into three types in \cite{LiShen}. The classification of projectively flat Finsler metrics with constant flag curvature were given in \cite{Li2} \cite{ShenProj0}.
Thus a natural problem is how to find more non-trivial projectively flat Finsler metrics which are not ($\alpha$, $\beta$)-metrics?

  A more general class named {\it general ($\alpha$, $\beta$)-metric} was first introduced by C. Yu and H. Zhu in \cite{Yu1} in the following form.
\be
F = \alpha \phi (b^2, \frac{\beta}{\alpha}),
\ee
where $\alpha$ is a Riemannian metric, $\beta$ is a 1-form, $b: = \| \beta_x\|_{\alpha}$ and $\phi = \phi(b^2, s)$ is a smooth function. It is easy to see that if $\phi_1 =0$, then $F$ is just a ($\alpha$, $\beta$)-metric. 
In \cite{Yu1}, C. Yu and H. Zhu found a class of projectively flat Finsler metrics in this class with the following three conditions: i) $\alpha = \sqrt{a_{ij}(x)y^i y^j}$ is projectively flat; (ii) $\beta$ is closed and conformal with respect to $\alpha$, i.e., the covariant derivatives of $\beta$ with respect to $\alpha$ is $b_{i|j} = c(x) a_{ij}$;
(iii)
$\phi = \phi(b^2, s)$ satisfies
$\phi_{22} = 2 (\phi_1 - s \phi_{12}).$

    In this paper, we only need the condition that $\alpha$ is projectively flat and determine all the  projectively flat  general ($\alpha$, $\beta$)-metrics.
Firstly, we get the following main theorem.
\begin{thm} \label{thm_main}
Let $F = \alpha \phi(b^2, \frac{\beta}{\alpha})$ be a non-Riemannian Finsler metric on an $n$-dimensional manifold with
$n \geq 3$. If $\alpha$ is a projectively flat Riemann metric and $\phi_1 \neq 0$, then
 $F$ is projectively flat if and only if one of the following conditions holds

(i) $\beta$ is parallel with respect to $\alpha$;

(ii) There are two scalar functions $c=c(b^2)$ and $k =k(x)$ such that
\be \label{thmcondi1}
[ c b^2 - (c-1)s^2 ]\phi_{22} = 2 b^2 (\phi_1 - s \phi_{12})
\ee
and
\be \label{thmcondi2}
b_{i|j} = k c (b^2 a_{ij} -b_i b_j) + k b_i b_j.
\ee
In this case,
\be G^i=^{\alpha}G^i + k \alpha\Big\{ (c-1) \frac{(b^2-s^2)\phi_2}{2\phi}  + \frac{b^2(2 s \phi_1 + \phi_2)}{2\phi}  \Big\}y^i.\ee
\end{thm}

Obviously, when $c=1$, it is just the case  studied in \cite{Yu1}.
By the above theorem, we can see that there are many choices of the function $c = c(b^2)$ and $k= k(x)$.
Actually, we find the general solutions of (\ref{thmcondi1}) and (\ref{thmcondi2}). Then we immediately obtain the following theorem.
\begin{thm} \label{thm_main2}
Let $F = \alpha \phi(b^2, \frac{\beta}{\alpha})$ be a non-Riemannian Finsler metric on an $n$-dimensional manifold with
$n \geq 3$ such that $\phi_1 \neq 0$. Suppose that $\alpha$ is a projectively flat Riemann metric with constant sectional  curvature $\kappa$, i.e.
\[ \alpha_{\kappa} = \frac{\sqrt{(1+\kappa|x|^2)|y|^2 - \kappa \langle x,y \rangle^2}}{1+\kappa|x|^2}. \]
If $\beta$ is not parallel with respect to $\alpha_{\kappa}$,
then
 $F$ is projectively flat if and only if there exists a scalar function $c = c(b^2)$ such that  $\phi$ and $\beta$ satisfy (\ref{thmgeneralsol1}) and (\ref{thmbeta}) respectively. That is
 \be \label{thmgeneralsol1}
\phi(b^2,s) = f(\mu + \nu s^2) - 2 \nu s \int^s_0 f'( \mu + \nu z^2 )dz + g(b^2) s, \ee
where $f=f(t)$ and $g=g(t)$ are two arbitrary $\mathcal{C}^{\infty}$ functions, and
\[  \mu = - \int c \nu d(b^2), \ \ \ \ \ \nu = -e^{\int \frac{c-1}{b^2} d(b^2) }. \]

\be \label{thmbeta}
 \beta =  e^{-\int \frac{c-1}{2 b^2} d(b^2)} \frac{ \varepsilon \langle x,y \rangle + (1 + \kappa |x|^2)\langle a,y \rangle - \kappa \langle a,x \rangle\langle x,y \rangle}{(1 + \kappa |x|^2)^{\frac{3}{2}}}.
 \ee
In this case,
\[ G^i=^{\alpha}G^i + e^{-\int \frac{c-1}{2 b^2} d(b^2)} \frac{\varepsilon -\kappa \langle a,x \rangle }{c b^2\sqrt{1 + \kappa |x|^2}}\alpha\Big\{ (c-1) \frac{(b^2-s^2)\phi_2}{2\phi}  + \frac{b^2(2 s \phi_1 + \phi_2)}{2\phi}  \Big\}y^i.\]
\end{thm}

Specially,  when  $\alpha$ is an Euclidean metric $|y|$ and $\beta = \langle x,y \rangle$, we get the spherically symmetric metric
\[ F = |y| \phi(|x|^2,\frac{\langle x,y \rangle}{|y|} ). \]
Obviously, it is a special class of general ($\alpha$, $\beta$)-metric. In this case, $\alpha=|y|$ is projectively flat and
$b_{i|j} = \delta_{ij}$. Thus by Theorem \ref{thm_main} and \ref{thm_main2} we get the following corollary which was first proved in
\cite{Zhou}.
\begin{cor}
Let $F = |y| \phi(|x|^2, \frac{\langle x,y \rangle}{|y|})$ be a non-Riemannian Finsler metric on an $n$-dimensional manifold with
$n \geq 3$. If $\phi_1 \neq 0$, then
 $F$ is projectively flat if and only if
\be \label{c=2condi1}
\phi_{22} = 2 (\phi_1 - s \phi_{12}),
\ee
i.e.
\[ \phi(b^2,s) = f(b^2 -  s^2) + 2 s \int^s_0 f'( b^2 - z^2 )dz + g(b^2) s,  \]
where $f=f(t)$ and $g=g(t)$ are two arbitrary $\mathcal{C}^{\infty}$ functions.
\end{cor}

\section{Preliminaries}

  In this section, we give some definitions and lemmas needed in this paper.
The geodesics of a  Finsler metric $F = F(x,y)$ on an open domain $\mathcal{U} \subset R^n$ are
determined by the following ODEs:
\[ \ddot{x} + 2 G^i (x, \dot{x}) = 0,\]
where $ G^i = G^i(x,y)$ are called {\it geodesic coefficients }
given by
\[ G^i = \frac{1}{4} g^{il} \Big\{  [F^2]_{x^m y^l} y^m -[F^2]_{x^l} \Big\}. \]

   A Finsler metric $F = F(x,y)$ is said to be {\it projectively flat } if its geodesic coefficients $G^i$ satisfy
\[ G^i = P y^i, \]
where $P = F_{x^k}y^k/(2F)$ is called the {\it projective factor} of $F$.
In this paper, we mainly consider the  general ($\alpha$, $\beta$)-metrics.

\begin{prop} {\rm (\cite{Yu1})}
Let $M$ be an $n$-dimensional manifold. $F = \alpha \phi(b^2, \frac{\beta}{\alpha})$
®´ is a Finsler metric on $M$ for any
Riemannian metric $\alpha$ and 1-form $\beta$ with $\|\beta\|_{\alpha} < b_o$ if and only if $\phi =\phi(b^2, s)$ is a positive $C^{\infty}$ function
satisfying
\be \label{positivecondi}
\phi - s \phi_2 > 0, \ \ \   \phi - s \phi_2 + (b^2 - s^2)\phi_{22} > 0
\ee
when $n \geq 3$ or
\[ \phi - s \phi_2 + (b^2 - s^2)\phi_{22} > 0\]
when n = 2, where $s$ and $b$ are arbitrary numbers with $|s|\leq b < b_o$.
\end{prop}

For simplicity, let
\[ r_{ij} : =\frac{1}{2}(b_{i|j}+b_{j|i}), \ \ \ r_j := b^i r_{ij}, \ \ \ s_{ij} : =\frac{1}{2}(b_{i|j}-b_{j|i}), \ \ \ s_j=b^i s_{ij}, \]
\[ r_{00} : =r_{ij}y^i y^j, \ \ \ r_{i0} := r_{ij}y^j, \ \ \ r_0 : = r_i y^i, \ \ \ s_{i0} : =s_{ij}y^j, \ \ \ s_0=s_i y^i. \]
In \cite{Yu1},
the geodesic coefficients $G^i$ of a general ($\alpha$, $\beta$)-metric $F = \alpha \phi(b^2, \frac{\beta}{\alpha})$ were given by
\be \label{G_upi}
\begin{split}
G^i =& ^{\alpha}G^i + \alpha Q s^i_{\ 0} + \big\{ \Theta (-2\alpha Q s_0 + r_{00}+ 2 \alpha^2 R r) + \alpha \Omega (r_0+s_0) \big\} \frac{y^i}{\alpha}\\
& + \big\{ \Psi (-2\alpha Q s_0 + r_{00}+ 2 \alpha^2 R r) + \alpha \Pi (r_0+s_0) \big\} b^i - \alpha^2 R (r^i + s^i),
\end{split}
\ee
where $^{\alpha}G^i$ are geodesic coefficients of $\alpha$,
\[ Q = \frac{\phi_2}{\phi - s\phi_2}, \ \ \ R= \frac{\phi_1}{\phi - s\phi_2},\]
\[ \Theta = \frac{(\phi - s\phi_2)\phi_2 - s\phi \phi_{22}}{2\phi[ \phi -s\phi_2 +(b^2-s^2)\phi_{22}]}, \ \ \
\Psi = \frac{\phi_22}{2[\phi -s\phi_2 +(b^2-s^2)\phi_{22}]}, \]
\[ \Pi = \frac{(\phi - s\phi_2)\phi_{12} -s \phi_1 \phi_{22}}{(\phi -s\phi_2)[\phi - s\phi_2 +(b^2-s^2)\phi_{22}]}, \ \ \
\Omega = \frac{2 \phi_1}{\phi} - \frac{s \phi + (b^2-s^2)\phi_2}{\phi}\Pi. \]

\section{Sufficient Conditions}
In this section, we  show the sufficient conditions for a general ($\alpha$, $\beta$)-metric to be   projectively flat. These conditions are also valid in dimension two.
\begin{lem} \label{lemsuffi}
Let $F = \alpha \phi(b^2, \frac{\beta}{\alpha})$ be a Finsler metric on an $n$-dimensional manifold with
$n \geq 2$. Suppose that $\beta$ satisfies

\be \label{lemscondi2}
b_{i|j} = k c (b^2 a_{ij} -b_i b_j) + k b_i b_j,
\ee
where $c=c(b^2)$ and $k =k(x)$ are two scalar functions.
Assume that $\phi = \phi(b^2,s)$ satisfies the following PDE:
\be \label{lemscondi1}
(c b^2 - (c-1)s^2)\phi_{22} = 2 b^2 (\phi_1 - s \phi_{12}).
\ee
Then the geodesic coefficients $G^i = G^i(x,y)$ of $F$ are given by
\[ G^i =^{\alpha}G^i + k  \alpha\Big\{(c-1) \frac{(b^2-s^2)\phi_2}{2\phi}  + \frac{b^2(2 s \phi_1 + \phi_2)}{2\phi}  \Big\}y^i. \]
\end{lem}
{\bf Proof}:
By the assumption, we have
\[ r_{00} = k c (b^2 \alpha^2 - \beta^2) + k \beta^2, \ \ \ r_0 = k b^2 \beta, \ \ \ r^i = k b^2 b^i, \]
\[ r = k b^4, \ \ \ s_0 = 0, \ \ \ s^i_{\ 0} =0. \]
Substituting them into (\ref{G_upi}) yields
\[
\begin{split}
G^i =& ^{\alpha}G^i  + k \alpha \{  \Theta [ c (b^2 - s^2) +  s^2 + 2 b^4 R ] + b^2 s  \Omega  \} y^i \\
& + k \alpha^2 \{  \Psi [ c (b^2 - s^2) +  s^2 + 2 b^4 R ] +  b^2 s \Pi - b^2 R \}b^i.
\end{split}
\]
By substituting the expression of $\Theta$, $\Psi$, $R$, $\Omega$ and $\Pi$ into the above equation we get
\[
\begin{split}
G^i  =& ^{\alpha}G^i  + k \alpha \Big\{ \frac{(b^2-s^2)(s \phi_{22}+\phi_2)}{2 \phi}(c-1) + \frac{ (s \phi_{22}+ 2 s^2 \phi_{12}+\phi_2) b^2}{2 \phi}   \\
&+  \frac{b^2 (2 s \phi -2 s^2\phi_2+ b^2 \phi_2 +b^2 s \phi_{22}-s^3 \phi_{22} ) }{\phi(\phi-s \phi_2+ b^2\phi_{22}-s^2 \phi_{22})}
[ \phi_1 - s \phi_{12} - \frac{c(b^2-s^2)+s^2}{2b^2}\phi_{22}] \Big\} y^i \\
& + k \alpha^2  \frac{ (c b^2 - (c-1)s^2)\phi_{22} -2 b^2 (\phi_1 - s \phi_{12}) }{2( \phi-s \phi_2+(b^2-s^2) \phi_{22} ) } b^i.
\end{split}
\]
Then by (\ref{lemscondi1}), we obtain
\[ G^i =^{\alpha}G^i + k \alpha\Big\{ \frac{(b^2-s^2)\phi_2}{2\phi} (c-1) + \frac{b^2(2 s \phi_1 + \phi_2)}{2\phi}  \Big\}y^i. \]
\Qed

{\bf Proof of the sufficiency of Theorem \ref{thm_main}}:
By the assumption, $\alpha$ is projectively flat, then $^{\alpha}G^i = ^{\alpha}P y^i$, where $^{\alpha}P$ is projective factor of $\alpha$.
Thus $F$ is a projectively flat Finsler metric by Lemma \ref{lemsuffi}.
\Qed

\section{Necessary Conditions}
The necessity of Theorem \ref{thm_main} can be obtained by the following lemma.
\begin{lem}
Let $F = \alpha \phi(b^2, \frac{\beta}{\alpha})$ be a non-Riemannian Finsler metric on an $n$-dimensional manifold with
$n \geq 3$ and $\phi_1 \neq 0$. Assume that $\alpha$ is projectively flat and $\beta$ is not parallel with respect to $\alpha$. If $F$ is projectively flat,
then there are two scalar functions $c=c(b^2)$ and $k =k(x)$ such that
\be \label{lemcondi1}
(c b^2 - (c-1)s^2)\phi_{22} = 2 b^2 (\phi_1 - s \phi_{12})
\ee
and
\be \label{lemcondi2}
b_{i|j} = k c (b^2 a_{ij} -b_i b_j) + k b_i b_j.
\ee
\end{lem}
{\bf Proof}: By the assumption, $\alpha$ is projectively flat, then $^{\alpha}G^i = ^{\alpha}P y^i$, where $^{\alpha}P$ is the projective factor of $\alpha$.  If $F$ is projectively flat, then by (\ref{G_upi}) there exists $\bar{P} = \bar{P}(x,y)$ such that
\be \label{lemP}
  \alpha Q s^i_{\ 0}+ \{\Psi (-2\alpha Q s_0 + r_{00}+ 2 \alpha^2 R r) + \alpha \Pi (r_0+s_0)\}b^i - \alpha^2 R (r^i + s^i) = \bar{P} y^i.
\ee
Contracting the above equation with  $y_i = a_{ij} y^j$, we have
\be
\bar{P} = \frac{s}{\alpha}\{ \Psi (-2\alpha Q s_0 + r_{00}+ 2 \alpha^2 R r) + \alpha \Pi (r_0+s_0)\} -  R (r_0 + s_0).
\ee
Substituting it back into (\ref{lemP}) yields
\be \label{lemPeq}
\begin{split}
&\big\{\Psi (-2\alpha Q s_0 + r_{00}+ 2 \alpha^2 R r) + \alpha \Pi (r_0+s_0)\big\}(b^i - \frac{s}{\alpha}y^i) \\
&+ \alpha Q s^i_{\ 0} -  R \big\{\alpha^2(r^i + s^i) + (r_0 + s_0)y^i \big\}  = 0.
\end{split}
\ee
Since the dependence of $\phi$ on $b^2$ and $s$ is unclear, it is difficult to solve the above equation directly. To overcome this problem, we choose a special coordinate system at a point as in  \cite{LiShen}.
Fix an arbitrary point $x_o \in \mathcal{U} \subset R^n$. Make a
change of coordinates: $(s,y^a) \rightarrow (y^i)$ by
\begin{eqnarray*}
y^1 = \frac{s}{\sqrt{b^2-s^2}}\bar{\alpha},\ \ \ y^a = y^a,\ \ \  a=2,...,n
\end{eqnarray*}
where $
\bar{\alpha}:= \sqrt{\sum^n_{a=2} (y^a)^2}.$
Then
\begin{eqnarray*}
\alpha = \frac{b}{\sqrt{b^2-s^2}}\bar{\alpha}, \ \ \ \beta= \frac{bs}{\sqrt{b^2-s^2}}\bar{\alpha}.
\end{eqnarray*}
and
\[ r_{00} = r_{11} \frac{s^2}{b^2-s^2}\bar{\alpha}^2 + 2 \bar{r}_{10}\frac{s}{\sqrt{b^2-s^2}}\bar{\alpha} + \bar{r}_{00},\]
\begin{eqnarray*}
r_0 = \frac{b s r_{11}}{\sqrt{b^2-s^2}}\bar{\alpha} + b \bar{r}_{10}, \ \ \  s_{0} = \bar{s}_0,
\end{eqnarray*}
\begin{eqnarray*}
s^1_{\ 0} = \bar{s}^1_{\ 0}, \ \ \  s^a_{\ 0} = \frac{s \bar{s}^a_{\ 1}}{\sqrt{b^2-s^2}}\bar{\alpha}+ \bar{s}^a_{\ 0},
\end{eqnarray*}
where $\bar{r}_{10} = r_{1a}y^a$, $\bar{s}_0 = s_a y^a$, $\bar{s}^1_{\ 0}=s_{1a} y^a$ and $\bar{s}^a_{\ 0} =s^a_{\ b} y^b$.
By a direct computation (\ref{lemPeq}) is equivalent to
\be     \label{lemEq1_1}
\begin{split}
&  (b^2-s^2) \Psi \bar{r}_{00} +  \Big\{ b^2  [ (b^2-s^2) (\Pi-2Q \Psi) + Q + s R ] \bar{s}_{10} +  [ (b^2-s^2)(b^2 \Pi + 2s \Psi) \\
&+ b^2 s R ] \bar{r}_{10} \Big\}\frac{\bar{\alpha}}{\sqrt{b^2-s^2}} +r_{11} (2 b^4 R \Psi + s^2 \Psi + b^2 s \Pi - b^2 R ) \bar{\alpha}^2 =0,
\end{split}
\ee
and
\be    \label{lemEq1_a}
\begin{split}
& s \Psi (b^2-s^2)y_a \bar{r}_{00}
+ \sqrt{b^2-s^2} \Big\{ b^2(-2 s Q \Psi - R + s \Pi) \bar{s}_{10} \\
&+ ( b^2 s \Pi + 2 s^2\Psi -b^2 R)  \bar{r}_{10} \Big\} y_a \bar{\alpha}   \\
&-\Big\{  b^2 \bar{s}_{a0} Q -y_a s r_{11} (2 b^4 R \Psi + s^2 \Psi + b^2 s \Pi - b^2 R  ) \Big\} \bar{\alpha}^2 \\
&
+ \Big\{ \frac{b^2}{\sqrt{b^2-s^2}} [ b^2  R (r_{1a}+ s_{1a}) + s Q s_{1a}] \Big\} \bar{\alpha}^3 =0. \end{split}
\ee
By $\bar{\alpha}= \sqrt{\sum^n_{a=2} (y^a)^2}$, (\ref{lemEq1_1}) is equivalent to
\be \label{lemequiv1}
(b^2-s^2) \Psi \bar{r}_{00} +r_{11} (2 b^4 R \Psi + s^2 \Psi + b^2 s \Pi - b^2 R ) \bar{\alpha}^2 =0
\ee
and
\be \label{lemequiv2}
 b^2  [ (b^2-s^2) (\Pi-2Q \Psi) + Q + s R ] \bar{s}_{10} +  [ (b^2-s^2)(b^2 \Pi + 2s \Psi)
+ b^2 s R ] \bar{r}_{10}  =0.
\ee
(\ref{lemEq1_a}) is equivalent to
\be \label{lemequiv3}
s \Psi (b^2-s^2)y_a \bar{r}_{00}-\Big\{  b^2 Q \bar{s}_{a0} -y_a s r_{11} (2 b^4 R \Psi + s^2 \Psi + b^2 s \Pi - b^2 R  ) \Big\} \bar{\alpha}^2 =0
\ee
and
\be \label{lemequiv4}
\begin{split}
& \sqrt{b^2-s^2} \Big\{ b^2(-2 s Q \Psi - R + s \Pi)\bar{s}_{10}
+ ( b^2 s \Pi + 2 s^2\Psi -b^2 R) \bar{r}_{10} \Big\} y_a\\
&+  \frac{b^2}{\sqrt{b^2-s^2}} [ b^2  R (r_{1a}+ s_{1a}) + s Q s_{1a}]  \bar{\alpha}^2 =0.
\end{split}
\ee

(\ref{lemequiv1})$\times s y_a - $(\ref{lemequiv3}) yields
\be
b^2 Q \bar{s}_{a0}\alpha^2 =0.
\ee
Then by the arbitrary choice of $x_{o}$,  $Q =0$ or $s_{ab}=0$. By the assumption that $F$ is non-Riemannian, we have $\phi_2 \neq 0$. Then
$Q \neq 0$. Thus $s_{ab} =0$.

Differentiating  (\ref{lemequiv2}) with respect to $y^a$ yields
\be \label{lemeq2_1}
 b^2  [ (b^2-s^2) (\Pi-2Q \Psi) + Q + s R ] s_{1a} +  [ (b^2-s^2)(b^2 \Pi + 2s \Psi)
+ b^2 s R ] r_{1a}  =0.
\ee

By (\ref{lemequiv4}) it can be seen that
\be \label{lemeq4_2}
(b^2 R + s Q) s_{1a}+ b^2  R r_{1a}  = 0
\ee
and
\be \label{lemeq4_1}
 b^2(-2 s Q \Psi - R + s \Pi) s_{1a}
+ ( b^2 s \Pi + 2 s^2\Psi -b^2 R) r_{1a} =0
\ee
because $\alpha^2$ is indivisible by $y_a$. We claim that $s_{1a}=0$.
If $R = 0$ at $x_{0}$, then by (\ref{lemeq4_2}) and $Q \neq 0$, we get $s_{1a} = 0$.
If $R \neq 0$ and $s_{1a}\neq 0$ at $x_{0}$, then (\ref{lemeq4_2}) becomes to
\[ b^2  + \frac{s Q}{R} + b^2   \frac{r_{1a}}{s_{1a}}  = 0.\]
 Differentiating the above equation with respect to $s$ yields
\[
\Big[ \frac{s \phi_2}{\phi_1}\Big]_2 = \frac{\phi_1\phi_2 + s \phi_1 \phi_{22} - s\phi_2 \phi_{12} }{\phi_1^2} =0.
\]
Then
\be \label{lemeqphi1}
 s \phi_1 \phi_{22} + \phi_2 (\phi_1 - s\phi_{12}) =0. \ee

It is easy to see that (\ref{lemeq4_2}) and (\ref{lemeq4_1}) both are linear equations of $s_{1a}$ and $r_{1a}$.
If $s_{1a}\neq 0$, then
\[
\begin{vmatrix}
\begin{matrix} (b^2 R + s Q) \end{matrix} &
\begin{matrix}b^2  R\end{matrix}
\\
\\
\begin{matrix} b^2(-2 s Q \Psi - R + s \Pi)\end{matrix} &
\begin{matrix}b^2 s \Pi + 2 s^2\Psi -b^2 R\end{matrix}
\end{vmatrix} =0.
\]
It is equivalent to
\be \label{lemeqphi2}
s (b^2 \phi_1 + s \phi_2) \phi_{22} + b^2 \phi_2 (-\phi_1+ s \phi_{12} ) =0.
\ee
(\ref{lemeqphi1})$\times b^2 +$(\ref{lemeqphi2}) yields
\be
s ( 2 b^2 \phi_1 + s \phi_2) \phi_{22} =0.
\ee
Then $\phi_{22} =0$ or $2 b^2 \phi_1 + s \phi_2=0$.
If $\phi_{22} =0$, then by (\ref{lemeqphi2}) we get
\[ \phi_2 (-\phi_1+ s \phi_{12} ) =0. \]
By the assumption that $F = \alpha \phi$ is non-Riemannian, we have $\phi_2 \neq 0.$
Then $\phi_1 =0$ when $s=0$. It is excluded. Thus
$\phi_{22} \neq 0$ and
\[ 2 b^2 \phi_1 + s \phi_2 =0.\]
Then $\phi_1 =0$ when $s=0$. It is excluded. Thus $s_{1a}=0$.

Substituting $s_{1a}=0$ into (\ref{lemeq4_2}) and (\ref{lemeq4_1}) yields
\be \label{s_10zero1}
 b^2 R r_{1a} =0
\ee
and
\be\label{s_10zero2}
 ( b^2 s \Pi + 2 s^2\Psi -b^2 R) r_{1a} =0.
\ee
If $R =0$, then $\phi_1 = 0$ and $\phi_{12} =0$. Then (\ref{s_10zero2}) becomes to
\[
 2 s^2\Psi r_{1a} =0.
\]
By the assumption that $\phi_{22} \neq 0$ which implies $\Psi\neq 0$, then we get
\[r_{1a} =0.\]

Differentiating (\ref{lemequiv1}) with respect to $y^a$ and $y^b$ yields
\[
(b^2-s^2) \Psi r_{ab} +r_{11} (2 b^4 R \Psi + s^2 \Psi + b^2 s \Pi - b^2 R ) \delta_{ab} =0.
\]
Set $r_{11} = k b^2$, ($k=k(x)$). Then there exists a number $c$ such that
\[ r_{ab} = c k b^2 \delta_{ab}. \]
Thus we obtain (\ref{lemcondi2}).

In this case, (\ref{lemequiv1}) becomes
\[ (b^2-s^2) \Psi c k +  (2 b^4 R \Psi + s^2 \Psi + b^2 s \Pi - b^2 R )k =0. \]
Substituting the expressions of $\Psi$, $\Pi$ and $R$ into the above equation yields
\be
k\frac{[c b^2 - (c-1)s^2]\phi_{22} - 2 b^2 (\phi_1 - s \phi_{12})}{2[\phi - s\phi_2 + (b^2-s^2)\phi_{22}]} =0.
\ee
If $k=0$, then $b_{i|j} =0$. It is contrary to the assumption. Thus we obtain (\ref{lemcondi2}).
\Qed

\section{General solutions of (\ref{thmcondi1})}\label{SectionC1}
To determine the projectively flat metrics in Theorem \ref{thm_main}, the efficient way is to solve  (\ref{thmcondi1}), (\ref{thmcondi2}). In this section,
we first give the general solutions of (\ref{thmcondi1}), then construct some special explicit solutions.
\begin{prop} \label{c=c(b^2)}
The general solutions of   (\ref{thmcondi1}) are given by
\be \label{generalsol1}
\phi(b^2,s) = f(\mu + \nu s^2) - 2 \nu s \int^s_0 f'( \mu + \nu z^2 )dz + g(b^2) s, \ee
where $f=f(t)$ and $g=g(t)$ are two arbitrary $\mathcal{C}^{\infty}$ functions,
\be \label{munu} \mu = - \int c \nu d(b^2), \ \ \ \ \ \nu = -e^{\int \frac{c-1}{b^2} d(b^2) }. \ee
\end{prop}
{\bf Proof}: Consider the following variable substitution,
\[ u = \mu + \nu s^2,\ \ \  v =s,   \]
where $\mu$ and $\nu$ are given by (\ref{munu})
By the implicit differentiation, we have
\begin{equation*}
\begin{split}
 0 = & \frac{\partial }{\partial v}\mu + s^2 \frac{\partial }{\partial v}\nu + 2s \nu \\
   = & - c \nu \frac{\partial (b^2)}{\partial v}  + s^2 \nu  \frac{c-1}{b^2} \frac{\partial (b^2)}{\partial v} + 2s. \nu
 \end{split}
\end{equation*}
Then
\[  \frac{\partial (b^2)}{\partial v} = \frac{2 s b^2}{ c b^2 - (c-1)s^2 }. \]
By the above equation and (\ref{thmcondi1}),
we get
\be \label{phizero}
\begin{split}
\frac{\partial}{\partial v} (\phi - s\phi_2) =& (\phi_1 - s \phi_{12}) \frac{\partial (b^2)}{\partial v} - s \phi_{22} \\
 = & (\phi_1 - s \phi_{12})  \frac{2 s b^2}{ c b^2 - (c-1)s^2 } - s \phi_{22} \\
 =& 0.
\end{split}
\ee
This means that
\be
\phi - s \phi_2 = f (\mu + \nu s^2)\label{phi2}
\ee
  is a function
of $\mu + \nu s^2$. Differentiating (\ref{phi2}) with respect to $s$ yields
\[ \phi_{22} = 2 \nu f'(\mu+\nu s^2) .\]
Thus
\[   \phi_2 = 2 \nu \int_0^s f'(\mu+\nu z^2) dz  + g (b^2).\]
Plugging it into (\ref{phi2})
gives (\ref{generalsol1}).
 \Qed

\bigskip

Let $c = \lambda = const.$, then (\ref{thmcondi1}) becomes
\be \label{c=lambda}
(\lambda b^2 - (\lambda -1)s^2)\phi_{22} = 2 b^2 (\phi_1 - s \phi_{12}).
\ee
By Proposition \ref{c=c(b^2)}, we immediately have
\begin{lem} \label{LemSpecial}
The general solutions of (\ref{c=lambda}) are given by
\[ \phi(b^2,s) = f(b^{2\lambda} - b^{2(\lambda-1)} s^2) + 2 b^{2(\lambda-1)} s \int^s_{0} f'(b^{2\lambda} - b^{2(\lambda-1)} z^2) dz + g(b^2) s,  \]
where $f=f(t)$ and $g=g(t)$ are two arbitrary $\mathcal{C}^{\infty}$ functions.
\end{lem}

In this case,
\[ \phi - s \phi_{2} = f(t), \ \ \ \ \phi - s \phi_{2} + (b^2-s^2)\phi_{22} = f(t) +  2 t f'(t), \]
where $t =b^{2(\lambda-1)}( b^{2} - s^2)$. Let $f(0) > 0$ and $f'(t)\geq 0$, then $\phi$ satisfies (\ref{positivecondi}). Thus $F = \alpha \phi$ is a Finsler metric.
Then we can construct many projectively flat general ($\alpha$, $\beta$)-metrics by choosing special $f$ and $g$.
The following are some special solutions.

(i) $f = 1$,
\[ \phi = 1  + g(b^2)s. \]
In this case, $F = \alpha \phi(b^2,s)$ is a Randers metric.

(ii) $f = \frac{1}{\sqrt{1-t}}$,
\[ \phi =\frac{\sqrt{1 - b^{2\lambda}+ b^{2(\lambda-1)} s^2}}{1-b^{2\lambda}}   + g(b^{2})s . \]
In this case, $F = \alpha \phi(b^2,s)$ is also a Randers metric.

(iii) $f=1+t$,
\[ \phi =b^{2 (\lambda-1)} s^2+ g(b^{2\lambda}) s+ 1+b^{2\lambda}. \]

(iv) $f=1+t^2$,
\[ \phi = -\frac{b^{4(\lambda-1)}}{3} s^4 + 2 b^{4\lambda -2} s^2+g(b^2) s + 1+b^{4\lambda}. \]

(v) $f = \ln(1+t)$,
\[ \phi = g(b^{2})s + \frac{2 b^{\lambda-1} \tanh^{-1} \Big[\frac{b^{\lambda-1} s}{ \sqrt{1 +b^{2\lambda}}  }\Big]} { \sqrt{1 +b^{2\lambda}} } s
+\ln(1+b^{2 \lambda}-b^{2 (\lambda-1)} s^2) \]
It is easy to see that (iii), (iv) and (v) are all new solutions.

\section{General solutions of (\ref{thmcondi2})} \label{SectionC2}
In this section we  solve  (\ref{thmcondi2}) in Theorem \ref{thm_main} when $\alpha$ is a projectively flat Riemannian metric. By Beltrami's theorem, projectively flat Riemannian metrics are equivalent to Riemannian metrics with constant sectional curvature. Let $\alpha_{\kappa}$ be a Riemannian metric with constant sectional curvature $\kappa$, then there is a local coordinate system  such that
\[ \alpha_{\kappa} = \frac{\sqrt{(1+\kappa|x|^2)|y|^2 - \kappa \langle x,y \rangle^2}}{1+\kappa|x|^2}. \]
Our method is to
 take a deformation of $\beta$ first, then find a conformal 1-form with respect to $\alpha_{\kappa}$. The following lemma can be obtained by a
 direct computation.
\begin{lem}
Let $\tilde{\beta} = \rho(b^2) \beta$, then
\[
\tilde{b}_{i|j} = \rho b_{i|j} + 2\rho' b_i(r_j + s_j). \]
\end{lem}
Let
$\rho = e^{\int \frac{c-1}{2 b^2} d(b^2)}$, then by  (\ref{thmcondi2}) and the above lemma we have
\[ \tilde{b}_{i|j} =c k b^2 \rho a_{ij}. \]
Then by the result in \cite{Xing}, we have
\[ \tilde{\beta} = \frac{ \varepsilon \langle x,y \rangle + (1 + \kappa |x|^2)\langle a,y \rangle - \kappa \langle a,x \rangle\langle x,y \rangle}{(1 + \kappa |x|^2)^{\frac{3}{2}}},
 \]
 \[ \tilde{b}_{i|j} = \frac{\varepsilon -\kappa \langle a,x \rangle }{\sqrt{1 + \kappa |x|^2}}a_{ij}, \]
where $\varepsilon$ is a constant number and $a$ is a constant vector.
When $c =0$, then $\varepsilon =0$ and $a = 0$. In this case, $\tilde{\beta} = \beta =0$.
Thus, when $c\neq 0$, we obtain
\begin{prop} \label{prop2}
If $\beta$ is a 1-form satisfies (\ref{thmcondi2}), where $c \neq 0$ and $\alpha_{\kappa} = \sqrt{a_{ij}y^i y^j}$ is a Riemannian metric with constant sectional flag curvature $\kappa$, then
\[ \beta =  e^{-\int \frac{c-1}{2 b^2} d(b^2)} \frac{ \varepsilon \langle x,y \rangle + (1 + \kappa |x|^2)\langle a,y \rangle - \kappa \langle a,x \rangle\langle x,y \rangle}{(1 + \kappa |x|^2)^{\frac{3}{2}}}. \]
In this case, the scalar function $k=k(x)$ in (\ref{thmcondi2}) is given by
\[ k(x) =  e^{-\int \frac{c-1}{2 b^2} d(b^2)} \frac{\varepsilon -\kappa \langle a,x \rangle }{c b^2\sqrt{1 + \kappa |x|^2}}.\]
\end{prop}

\vskip 10mm
\noindent
Benling Li\\
Department of Mathematics\\
Ningbo University\\
Ningbo, Zhejiang Province 315211\\
P.R. China\\
libenling@nbu.edu.cn

\vskip 10mm
\noindent
Zhongmin Shen\\
Department of Mathematical Sciences\\
Indiana University Purdue University Indianapolis\\
402 N. Blackford Street\\
Indianapolis, IN 46202-3216, USA\\
zshen@math.iupui.edu


\begin{thebibliography}{00}
\bibitem{Be2} L. Berwald,  {\it $\ddot{U}$ber die n-dimensionalen  Geometrien konstanter Kr$\ddot{u}mmung$, in denen die Geraden die $k\ddot{u}rzesten$ sind,} Math. Z.\textbf{ 30}(1929), 449-469.
\bibitem{Bryant} R. Bryant, {\it Finsler structures on the 2-sphere satisfying $K = 1$}, Finsler Geometry,
Contemporary Mathematics \textbf{196}, Amer. Math. Soc., Providence,
RI, 1996, 27-42. MR \textbf{97e}:53128.
\bibitem{Bry2} R. Bryant, {\it Projectively flat Finsler 2-spheres of constant curvature}, Selecta Math.,
New Series, \textbf{3}(1997), 161-204. MR \textbf{98i}:53101.
\bibitem{Chen} B. Chen, Z. Shen and L. Zhao {\it On a class of Ricci-flat Finsler metrics in Finsler geometry}, J. Geom. Phys. \textbf{70} (2013), 30¨C38.
\bibitem{Cheng} X. Cheng, Z. Shen and Y. Zhou, {\it On a class of locally dually flat Finsler Metrics}, Int. J. Math., {\bf 21} (2010), 1531-1543.
\bibitem{Li1} B. Li, {\it On dually flat Finsler metrics}, Differ. Geom. Appl. \textbf{31}(2013), 718-724.
\bibitem{Li2} B. Li, {\it On the classification of projectively flat Finsler metrics with constant flag curvature}, Adv. Math. \textbf{257}(2014), 266-284.

\bibitem{LiShen} B. Li and Z. Shen, {\it On a class of projectively flat Finsler metrics with constant flag curvature}, Int. J. Math., Vol.18, No. 7 (2007), 749-760.
\bibitem{Huang} L. Huang and X. Mo, {\it On some explicit constructions of dually flat Finsler metrics}
J. Math. Anal. Appl., Vol. 405, No. 2 (2013), 565-573.

\bibitem{Sevim} E. S. Pending Sevim, Z. Shen and L. Zhao {\it Some Ricci-flat Finsler metrics}, Publ. Math. Debrecen \textbf{83}(4) (2013), 617¨C623.
\bibitem{ShenProj0} Z.Shen, {\it Projectively flat Finsler metrics of constant flag curvature}, Trans. Amer. Math. Soc. \textbf{355}(4)(2003), 1713¨C11728.
\bibitem{ShenProj1} Z. Shen, {\it On Projectively flat $(\alpha,\beta)$-metrics}, Can. Math.
Bull., {\bf 52}(1)(2009), 132-144.
\bibitem{Xia} Q. Xia, {\it On locally dually flat ($\alpha$, $\beta$)-metrics}, Differ. Geom. Appl., {\bf 29} (2011), 233-243.
\bibitem{Xing} H. Xing and Z. Shen, {\it On Randers Metrics with Isotropic S-Curvature}, Acta Math. Sinica, {\bf 24}(5) (2008), 789-796.
\bibitem{Yu1} C. Yu and H. Zhu, {\it On a new class of Finsler metrics}, Differ. Geom. Appl. {\bf 29}(2011), 244-254.
\bibitem{Zhou} L. Shou, {\it Spherically symetric Finsler metrics in $R^n$}, Publ. Math. Debrecen.  {\bf 80} 1-2 (4) (2012), 67-77.



\end{thebibliography}
\end{document}